\documentclass{amsart}
\usepackage{amsmath}
\usepackage{array}
\usepackage{amsfonts}
\usepackage{latexsym}
\usepackage{epsfig}
\newtheorem{theorem}{Theorem}[section]
\newtheorem{prop}[theorem]{Proposition}

\newtheorem{lemma}[theorem]{Lemma}

\newcommand\beq{\begin{equation}}
\newcommand\eeq{\end{equation}}
\newcommand\bce{\begin{center}}
\newcommand\ece{\end{center}}
\newcommand\bea{\begin{eqnarray}}
\newcommand\eea{\end{eqnarray}}
\newcommand\ben{\begin{enumerate}}
\newcommand\een{\end{enumerate}}
\newcommand\bsb{\begin{Sb}}
\newcommand\esb{\end{Sb}}
\newcommand\us{\underset}

\newcommand\wt{\widetilde}

\newcommand\nn{\nonumber}

\newcommand\ms{\medskip}

\newcommand\brr{\begin{array}}
\newcommand\err{\end{array}}
\newcommand\bt{\begin{tabular}}
\newcommand\et{\end{tabular}}

\renewcommand\S{{\mathcal S}}
\newcommand\D{{\mathcal D}}
\newcommand\rs{\Phi_\llcorner}
\newcommand\bjs{\Psi_\llcorner}
\newcommand\kra{\Phi_\urcorner}
\newcommand\bijn{\Psi_\urcorner}
\newcommand\krat{\varphi_K}

\newenvironment{abstrac}{%
         \small
        \begin{center}%
          {\bfseries {Abstract}\vspace{-.5em}}%
                 \end{center}%
        \quotation}

\title{Fixed points and excedances in restricted permutations}
\author{Sergi Elizalde}
\address{Department
of Mathematics, Massachusetts Institute of Technology,  77
Massachusetts Avenue, Cambridge MA
02139.}\email{sergi@math.mit.edu}
\date{November 2002}
\begin{document}
 \maketitle
 \thispagestyle{empty}

\begin{abstrac}
\medskip
In this paper we prove that among the permutations of length $n$
with $i$ fixed points and $j$ excedances, the number of
321-avoiding ones equals the number of 132-avoiding ones, for all
given $i,j\le n$.
 We use a new technique involving diagonals of non-rational
generating functions.

This theorem generalizes a recent result of Robertson, Saracino
and Zeilberger, for which we also give another, more direct proof.
\end{abstrac}

\section{Introduction}

Let $n$, $m$ be two positive integers with $m\le n$, and let
$\pi=\pi_1\pi_2\cdots\pi_n\in\S_n$ and
$\sigma=\sigma_1\sigma_2\cdots\sigma_m\in\S_m$. We say that $\pi$
\emph{contains} $\sigma$ if there exist indices
$i_1<i_2<\ldots<i_m$ such that $\pi_{i_1}\pi_{i_2}\cdots\pi_{i_m}$
is in the same relative order as $\sigma_1\sigma_2\cdots\sigma_m$.
If $\pi$ does not contain $\sigma$, we say that $\pi$ is
\emph{$\sigma$-avoiding}. For example, if $\sigma=132$, then
$\pi=24531$ contains $\sigma$, because $\pi_1\pi_3\pi_4=253$.
However, $\pi=42351$ is $\sigma$-avoiding.

We say that $i$ is a \emph{fixed point} of a permutation $\pi$ if
$\pi_i=i$, and that $i$ is an \emph{excedance} of $\pi$ if
$\pi_i>i$. Denote by $\mathrm{fp}(\pi)$ and $\mathrm{exc}(\pi)$
the number of fixed points and the number of excedances of $\pi$
respectively. Denote by $\S_n(\sigma)$ the set of
$\sigma$-avoiding permutations in $\S_n$. We are interested in the
distribution of the number of fixed points and excedances among
the permutations in $\S_n(\sigma)$.

For the case of patterns of length 3, it is known (\cite{Knu73})
that regardless of the pattern $\sigma\in\S_3$,
$|\S_n(\sigma)|=C_n=\frac{1}{n+1}{2n\choose n}$, the $n$-th
Catalan number. Bijective proofs of this fact are given in
\cite{Kra01,Ric88,SS85,Wes95}.

In the recent paper \cite{RSZ02}, pattern-avoiding permutations
are studied with respect to the number of fixed points, and an
interesting refinement is presented. It is shown that given $i\le
n$, the number of 321-avoiding permutations of length $n$ with $i$
fixed points equals the number of 132-avoiding permutations of
length $n$ with $i$ fixed points.

In this paper we prove a further refinement of this result, namely
that it still holds when we fix not only the number of fixed
points but also the number of excedances. In other words, the
bivariate distribution of fixed points and excedances is the same
in both 321-avoiding permutations and in 132-avoiding
permutations.

One of the key points in the proof is to use bijections between
pattern avoiding permutations and Dyck paths. Recall that a
\emph{Dyck path} of length $2n$ is a lattice path in
$\mathbb{Z}^2$ between $(0,0)$ and $(2n,0)$ consisting of up-steps
$(1,1)$ and down-steps $(1,-1)$ which never goes below the
$x$-axis. Sometimes it will be convenient to encode each up-step
by a letter $u$ and each down-step by $d$, obtaining an encoding
of the Dyck path as a \emph{Dyck word}. We shall denote by $\D_n$
the set of Dyck paths of length $2n$, and by
$\D=\bigcup_{n\geq0}\D_n$ the class of all Dyck paths. It is
well-known that $|\D_n|=C_n$. If $D\in\D_n$, we will write $|D|=n$
to indicate the semilength of $D$. The generating function
(\emph{GF} for short) that enumerates Dyck paths according to
their semilength is $\sum_{D\in\D}{t^{|D|}}=\sum_{n\geq0}{C_n t^n}
=\frac{1-\sqrt{1-4t}}{2t}$, which we denote by $C(t)$.

\section{Statement of the main theorem}

Here is the main result of this paper.
\begin{theorem}\label{th:main}
For any $0\le i,j\le n$, \bea\nn|\{\pi\in\S_n(321):\ \mathrm{fp}(\pi)=i,\
\mathrm{exc}(\pi)=j\}|=|\{\pi\in\S_n(132):\ \mathrm{fp}(\pi)=i,\ \mathrm{exc}(\pi)=j\}|.\eea

\ms Equivalently,
$$\sum_{\pi\in\S_n(321)}{x^{\mathrm{fp}(\pi)}q^{\mathrm{exc}(\pi)}}=\sum_{\pi\in\S_n(132)}{x^{\mathrm{fp}(\pi)}q^{\mathrm{exc}(\pi)}}.$$
\end{theorem}

The proof of this theorem is done in two parts. First, in section
\ref{sec:321} we find the GF for the number of fixed points and
excedances in 321-avoiding permutations. Then, in section
\ref{sec:132}, we show that this GF also counts the number of
fixed points and excedances in 132-avoiding permutations. To do
the latter, we introduce an extra variable marking a new parameter
in the GF. Then, using combinatorial properties, we deduce an
identity that determines this GF. Finally, we conjecture an
expression for it and check that our expression satisfies the
identity, hence it is the correct GF.

\section{Counting 321-avoiding permutations according to fixed points and
excedances}\label{sec:321}

The goal of this section is to find an expression for the GF
$$F_{321}(x,q,t):=\sum_{n\geq0}\sum_{\pi\in\S_n(321)}x^{\mathrm{fp}(\pi)}q^{\mathrm{exc}(\pi)}t^n.$$
Instead of counting fixed points and excedances directly in
321-avoiding permutations, we define the following bijection $\rs$
between $\S_n(321)$ and $\D_n$, suggested by Richard Stanley.

\begin{figure}[hbt]
\epsfig{file=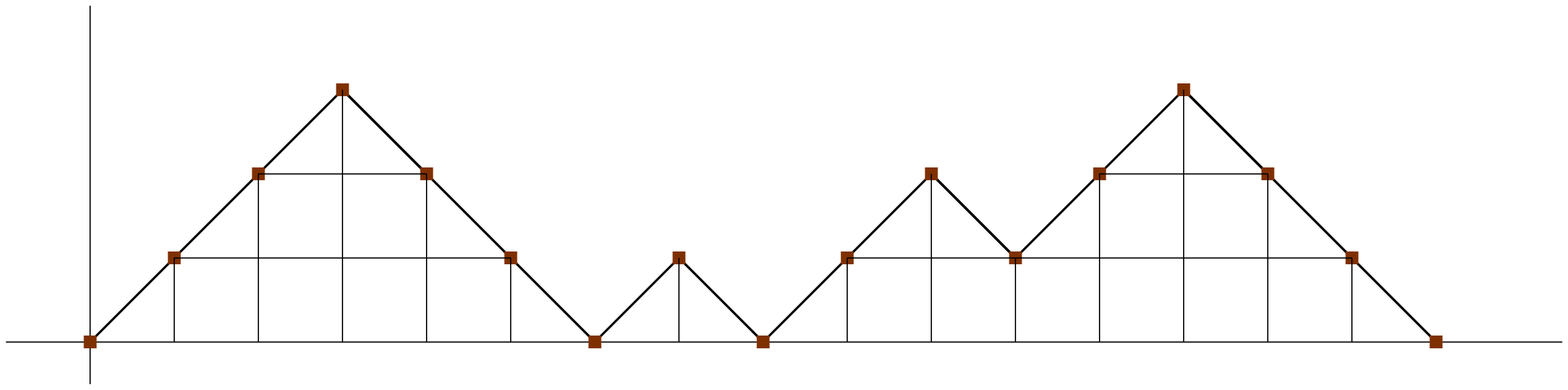,width=4.5in} \caption{\label{fig:d1} The Dyck
path $\rs(23147586)$.}
\end{figure}

Given $\pi=\pi_1\pi_2\cdots\pi_n\in\S_n(321)$, let
$a_i=\mathrm{max}\{j:\{1,2,\ldots,j\}\subseteq\{\pi_1,\pi_2,\ldots,\pi_i\}\}$
($j$ can be 0), for each $0\le i\le n$. Now build the Dyck path
$\rs(\pi)$ by adjoining, for each $i$ from 1 to $n$, one up-step
followed by $\mathrm{max}\{a_i-\pi_i+1,0\}$ down-steps. For
example, the Dyck path corresponding to $\pi=23147586$ is given in
Figure~\ref{fig:d1}.

There is an alternative way to define this bijection. A
\emph{right-to-left minimum} of $\pi$ is an element $\pi_i$ such
that $\pi_i<\pi_j$ for all $j>i$. Let
$\pi_{i_1},\pi_{i_2},\ldots,\pi_{i_k}$ be the right-to-left minima
of $\pi$, from left to right. For example, the right-to-left
minima of $23147586$ are $1,4,5,6$. Then, $\rs(\pi)$ is precisely
the path that starts with $i_1$ up-steps, then has, for each $j$
from 2 to $k$, $\pi_{i_j}-\pi_{i_{j-1}}$ down-steps followed by
$i_j-i_{j-1}$ up-steps, and finally ends with $n+1-\pi_{i_k}$
down-steps.

\begin{figure}[hbt]
\epsfig{file=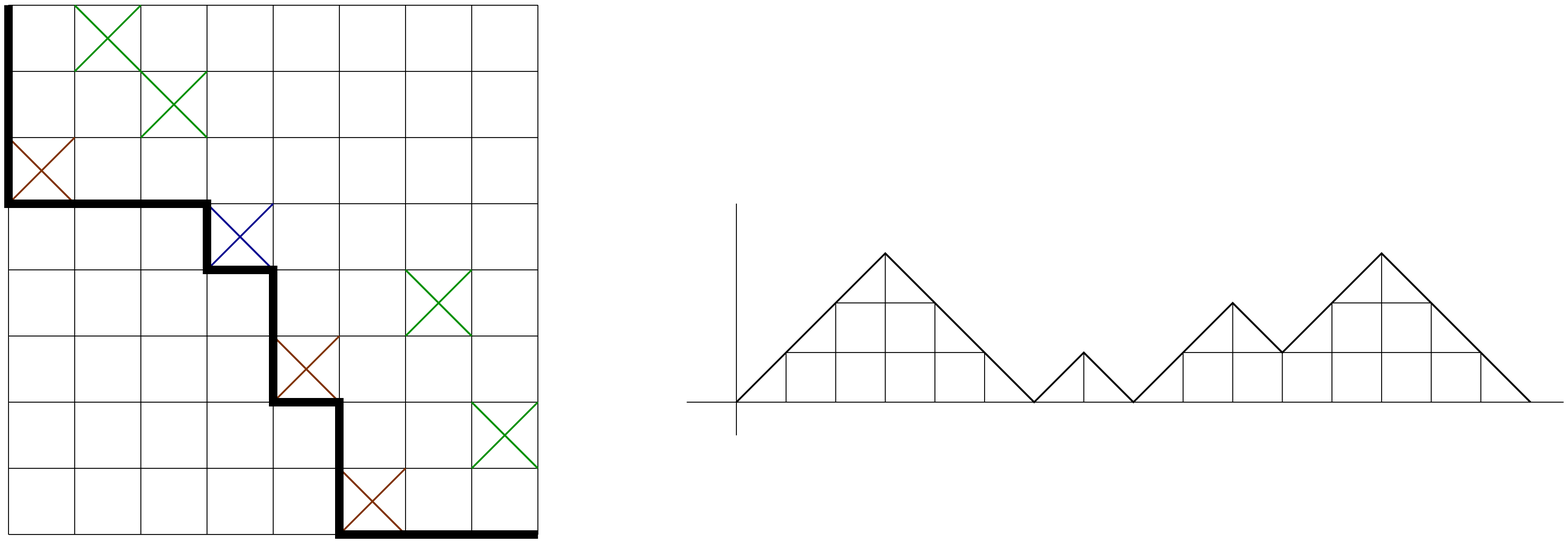,height=3.7cm} \caption{\label{fig:bij_rs}
The bijection $\rs$.}
\end{figure}

An easy way to picture this construction is to represent $\pi$ as
an $n\times n$ array with a cross on the squares $(i,\pi_i)$. It
is known that a permutation is 321-avoiding if and only if both
the subsequence determined by its excedances and the one
determined by the remaining elements are increasing. In this array
representation, excedances correspond to crosses strictly to the
right of the main diagonal. Note that the rest of the crosses are
precisely the right-to-left minima. Consider the path with
\emph{down} and \emph{right} steps along the edges of the squares
that goes from the upper-left corner to the lower-right corner of
the array leaving all the crosses to the right and remaining
always as close to the main diagonal as possible. Then $\rs(\pi)$
can be obtained from this path just by reading an up-step every
time the path moves down, and a down-step every time the path
moves to the right. Figure~\ref{fig:bij_rs} shows a picture of
this bijection, again for $\pi=23147586$.

Recall that a \emph{peak} of a Dyck path $D\in\D$ is an up-step
followed by a down-step (i.e., an occurrence of $ud$ in the
associated Dyck word). A \emph{hill} is a peak at height 1, where
the height is the $y$-coordinate of the top of the peak. Denote by
$h(D)$ the number of hills of $D$. A \emph{double rise} of a Dyck
path is an up-step followed by another up-step ($uu$ when seen as
a word). Denote by $\mathrm{dr}(D)$ the number of double rises of $D$.

It can easily be checked that $\rs$ has the property that
$\mathrm{fp}(\pi)=h(\rs(\pi))$ and $\mathrm{exc}(\pi)=\mathrm{dr}(\rs(\pi))$. Therefore,
counting 321-avoiding permutations according to the number fixed
points and excedances is equivalent to counting Dyck paths
according to the number of hills and double rises. More precisely,
$$F_{321}(x,q,t)=\sum_{D\in\D}{x^{h(D)}q^{\mathrm{dr}(D)}}t^{|D|}.$$

We can give an equation for $F_{321}$ using the symbolic method
described in \cite{FlSe98} and \cite{SeFl96}. A recursive
definition for the class $\D$ is given by the fact that every
non-empty Dyck path $D$ can be decomposed in a unique way as
$D=uAdB$, where $A,B\in\D$. Clearly, $h(D)=h(B)+1$ and
$\mathrm{dr}(D)=\mathrm{dr}(B)$ if $A$ is empty, and $h(D)=h(B)$ and
$\mathrm{dr}(D)=\mathrm{dr}(A)+\mathrm{dr}(B)+1$ otherwise. Hence, we obtain the following
equation for $F_{321}$: \bea
F_{321}(x,q,t)=1+t(x+q(F_{321}(1,q,t)-1))F_{321}(x,q,t).
\label{F321x1} \eea Substituting first $x=1$, we obtain that
$F_{321}(1,q,t)=\frac
{1+t(q-1)-\sqrt{1-2t(1+q)+t^2(1-q)^2}}{2qt}$. Now, solving
(\ref{F321x1}) for $F_{321}(x,q,t)$ gives \bea F_{321}(x,q,t)=
\frac{2}{1+t(1+q-2x)+\sqrt{1-2t(1+q)+t^2(1-q)^2}}. \label{F321}
\eea

To conclude this section, we want to remark that applying this
method one can also obtain the GF that enumerates fixed points,
excedances and descents in 321-avoiding permutations. It can be
seen that the number of descents of a permutation $\pi$ (i.e.,
indices $i$ for which $\pi_i>\pi_{i+1}$), denoted $\mathrm{des}(\pi)$,
equals the number of occurrences of $uud$ in the Dyck word of
$\rs(\pi)$. Using the same decomposition as before, we conclude
that
$$\sum_{n\geq0}\sum_{\pi\in\S_n(321)}x^{\mathrm{fp}(\pi)}q^{\mathrm{exc}(\pi)}p^{\mathrm{des}(\pi)}t^n=
\frac{2}{1+t(1+q-2x)+\sqrt{1-2t(1+q)+t^2((1+q)^2-4qp)}}.$$

\section{Counting 132-avoiding permutations according to fixed points and
excedances}\label{sec:132}

Analogously to the previous section, we define
$$F_{132}(x,q,t):=\sum_{n\geq0}\sum_{\pi\in\S_n(132)} x^{\mathrm{fp}(\pi)}q^{\mathrm{exc}(\pi)}t^n.$$
To prove Theorem~\ref{th:main} we have to show that
$F_{321}(x,q,t)=F_{132}(x,q,t)$.

Instead of enumerating fixed points and excedances directly in
132-avoiding permutations, we use a bijection between $\S_n(132)$
and $\D_n$, and then look at what are the statistics
in Dyck paths that correspond to $\mathrm{fp}$ and $exc$ after the bijection.

For any $D\in\D$, we define a \emph{tunnel} of $D$ to be a
horizontal segment between two lattice points
of $D$ that intersects $D$ only in these
two points, and stays always below $D$. Tunnels are in
obvious one-to-one correspondence with decompositions of the Dyck word $D=AuBdC$,
where $B\in\D$ (no restrictions on $A$ and $C$). In the decomposition, the tunnel is the segment
that goes from the beginning of $u$ to the end of $d$.
If $D\in\D_n$, then $D$ has exactly $n$ tunnels, since such a
decomposition can be given for each up-step of $D$.

A tunnel of $D\in\D_n$ is called a \emph{centered tunnel} if the
$x$-coordinate of its midpoint (as a segment) is $n$, that is, the
tunnel is centered with respect to the vertical line through the
middle of $D$. In terms of the decomposition $D=AuBdC$, this is
equivalent to saying that $A$ and $C$ have the same length. Denote
by $\mathrm{CT}(D)$ the set of centered tunnels of $D$, and let
$\mathrm{ct}(D)=|\mathrm{CT}(D)|$.

\begin{figure}[hbt]
\epsfig{file=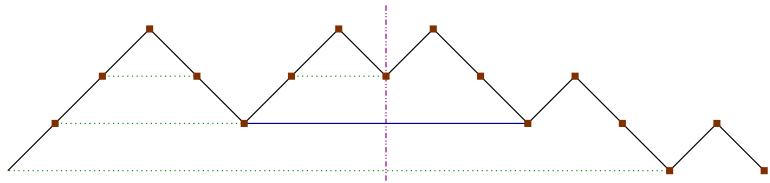,width=4in} \caption{\label{fig:ctlt}
Centered and left tunnels.}
\end{figure}

A tunnel of $D\in\D_n$ is called a \emph{left tunnel} if the
$x$-coordinate of its midpoint is strictly less than $n$, that is,
the midpoint of the tunnel is to the left of the vertical line
through the middle of $D$. In terms of the decomposition
$D=AuBdC$, this is equivalent to saying that the length of $A$ is strictly smaller
than the length of $C$. Denote by $\mathrm{lt}(D)$ the number of left tunnels of $D$.
In Figure~\ref{fig:ctlt}, there is one centered tunnel drawn with a solid line,
and four left tunnels drawn with dotted lines.

\begin{figure}[hbt]
\epsfig{file=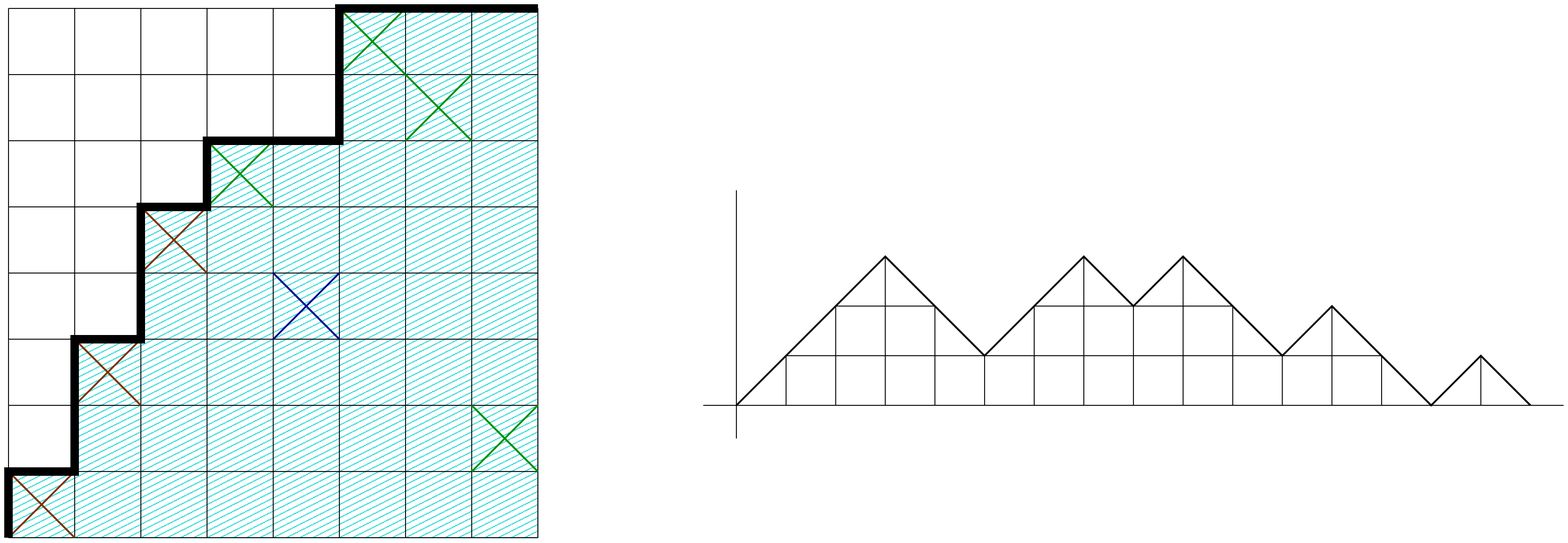,height=3.7cm} \caption{\label{fig:bij2}
The bijection $\krat$.}
\end{figure}

We will use the bijection between $\S_n(132)$ and $\D_n$ given
by Krattenthaler in \cite{Kra01}. We denote it by $\krat$. For
$\pi=\pi_1\pi_2\cdots\pi_n\in\S_n(132)$, $\krat(\pi)$ is obtained
by reading $\pi$ from left to right and adjoining for each $\pi_j$
as many up-steps as necessary followed by a down-step from height
$h_j+1$ to height $h_j$, where $h_j$ is the number of elements in
$\pi_{j+1}\cdots\pi_n$ which are larger than $\pi_j$. As pointed
out by Reifegerste in \cite{Rei02}, this path is closely related
to the diagram of $\pi$ obtained from the $n\times n$ array
representation of $\pi$ by shading, for each cross, the cell containing it
and the squares that are due south and due east of it.
The diagram, defined as the region that remains unshaded,
is determined by the path with \emph{left} and \emph{down} steps
that goes from the upper-right corner to the lower-left corner,
leaving all the crosses to the right, and staying always as close
to the diagonal connecting these two corners as possible. If we go along this path
reading an up-step every time it goes left and a down-step every time it goes down,
we get $\krat(\pi)$. Figure~\ref{fig:bij2} shows an example when
$\pi=67435281$.

What is interesting of this bijection for our purposes is that it
has the property that it maps fixed points to centered tunnels,
and excedances to left tunnels. This can be seen using the diagram
representation. There is an easy way to recover a permutation
$\pi\in\S_n(321)$ from its diagram: row by row, put a cross in the
leftmost shaded square such that there is exactly one dot in each
column. Now, instead of looking directly at $\krat(\pi)$, consider
the path from the upper-right corner to the lower-left corner of
the array of $\pi$. To each cross we can associate a tunnel in a
natural way. Indeed, each cross produces a decomposition
$\krat(\pi)=AuBdC$ where $B$ corresponds to the part of the path
above and to the left of the cross. Here $u$ corresponds to the
horizontal step directly above the cross, and $d$ to the vertical
step directly to the left of the cross. Thus, fixed points, which
correspond to crosses on the main diagonal, give centered tunnels,
and excedances, which are crosses to the right of the main
diagonal, give left tunnels. This means that
$\mathrm{fp}(\pi)=\mathrm{ct}(\krat(\pi))$ and
$\mathrm{exc}(\pi)=\mathrm{lt}(\krat(\pi))$. So, our problem is
equivalent to counting Dyck paths according to centered and left
tunnels, and the function we want to find becomes
$$F_{132}(x,q,t)=\sum_{D\in\D}{x^{\mathrm{ct}(D)}q^{\mathrm{lt}(D)}}t^{|D|}.$$

The decomposition of $\D$ that we used to enumerate hills and
double rises no longer works here. Indeed, if we write
$D=uAdB$ with $A,B\in\D$, then $\mathrm{ct}(A)$ and $\mathrm{ct}(B)$ do not give
information about $\mathrm{ct}(D)$. However, to count only centered
tunnels, we can use another decomposition.

Now we show how to obtain an expression for $F_{132}(x,1,t)$. We
consider Dyck paths with \emph{marked} centered tunnels. That is,
we count pairs $(D,S)$ where $D\in\D$ and $S\subseteq \mathrm{CT}(D)$. Each
such pair is given weight $(x-1)^{|S|} t^{|D|}$, so that for a fixed
$D$, the sum of weights of all pairs $(D,S)$ will be
$\sum_{S\subseteq \mathrm{CT}(D)}{(x-1)^{|S|} t^{|D|}}=((x-1)+1)^{|\mathrm{CT}(D)|}
t^{|D|}=x^{\mathrm{ct}(D)} t^{|D|}$, which is precisely the weight that $D$
has in $F_{132}(x,1,t)$.

\begin{figure}[hbt]
\epsfig{file=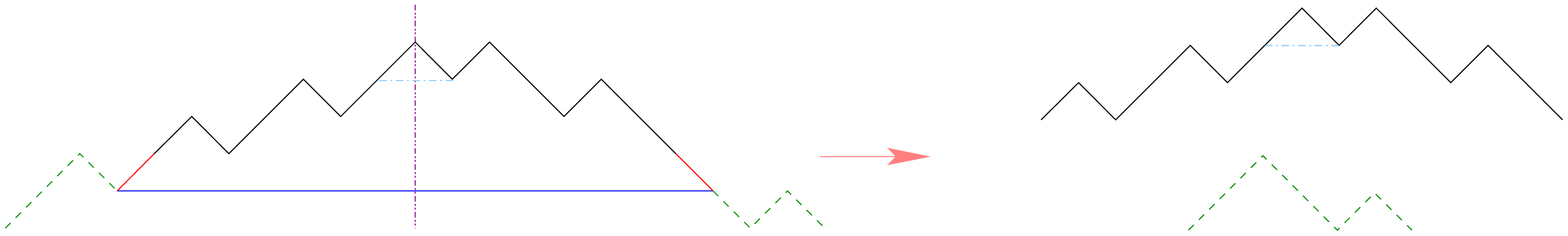,width=5.8in} \caption{\label{fig:countct}
Decomposing Dyck paths with marked centered tunnels.}
\end{figure}

Dyck paths with no marked tunnels (i.e., pairs $(D,\emptyset)$)
are enumerated by $C(t)=\frac{1-\sqrt{1-4t}}{2t}$, the GF for
Catalan numbers. On the other hand, for an arbitrary Dyck path
$D$ with some centered tunnel marked (i.e., a pair $(D,S)$ with
$S\neq\emptyset$), we can consider the decomposition given by the
longest marked tunnel, say $D=AuBdC$. Then, $AC$ (seen as the
concatenation of Dyck words) gives an arbitrary Dyck path with no
marked centered tunnels, and $B$ is an arbitrary Dyck path where
some centered tunnels may be marked (Figure~\ref{fig:countct}).
This decomposition translates into the following equation for GFs:
$$F_{132}(x,1,t)=C(t)+(x-1)t C(t) F_{132}(x,1,t).$$
Solving it, we obtain
$$F_{132}(x,1,t)=\frac{2}{1+2t(1-x)+\sqrt{1-4t}},$$ which is
precisely the expression that we had for $F_{321}(x,1,t)$ in
(\ref{F321}). This gives a new and perhaps simpler proof of the
main result in \cite{RSZ02}, namely that
$|\{\pi\in\S_n(321):\mathrm{fp}(\pi)=i\}|=|\{\pi\in\S_n(132):\mathrm{fp}(\pi)=i\}|$
for all $i\le n$.

\ms

To enumerate left tunnels we will need a different approach. The
first step is to generalize the concepts of centered and left
tunnels, allowing the vertical line that we use as a reference to
be shifted from the center of the Dyck path. For $D\in\D$ and
$r\in\mathbb{Z}$, let $\mathrm{ct}_r(D)$ be the number of tunnels of $D$
whose midpoint lies on the vertical line $x=n-r$ (we call this the
\emph{reference line}). Similarly, let $\mathrm{lt}_r(D)$ be the number of
tunnels of $D$ whose midpoint lies on the half-plane $x<n-r$.
Notice that by definition, $\mathrm{ct}_0$ and $\mathrm{lt}_0$ are respectively
the statistics $ct$ and $\mathrm{lt}$ defined previously.

We also add a new variable $v$ to $F_{132}$ which marks the
distance from the reference line to the actual middle of the path.
Define
$$G(x,q,t,v):=\sum_{n,r\geq0}\sum_{D\in\D_n}x^{\mathrm{ct}_r(D)}q^{\mathrm{lt}_r(D)}v^r t^n.$$

Our next goal is to find an equation that determines $G(x,q,t,v)$.
The idea is to use again the decomposition of a Dyck path as
$D=uAdB$, where $A,B\in\D$. The difference is that now the GFs involve sums not only over Dyck
paths but also over the possible positions of the reference line.

\begin{figure}[hbt]
\epsfig{file=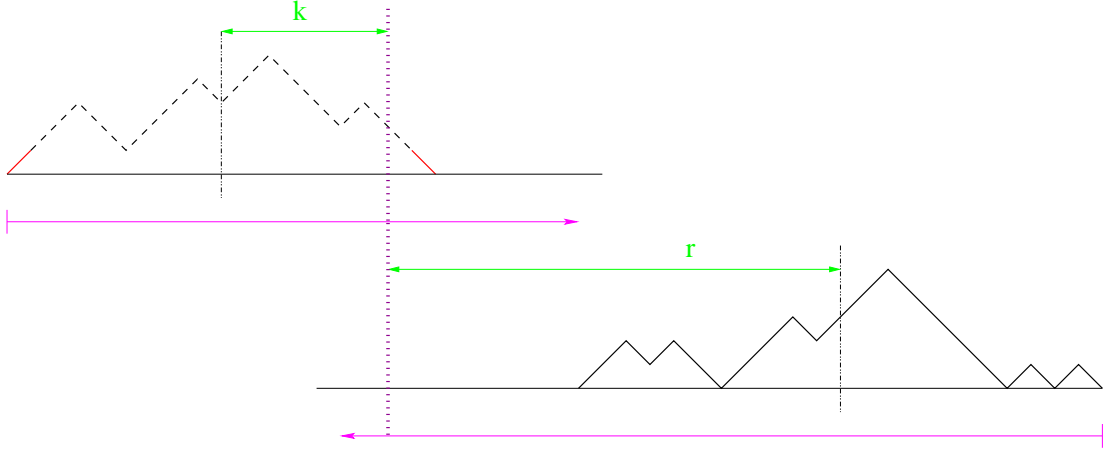,height=6cm}
 \caption{\label{fig:prod} $H_1$ and $H_2$.}
\end{figure}

Let \bea H_1(x,q,t,v):=
\us{k\ge-n}{\sum_{n\ge1}}\sum_{A\in\D_{n-1}}x^{\mathrm{ct}_{-k}(uAd)}q^{\mathrm{ct}_{-k}(uAd)}v^k
t^n\label{h1def}\eea be the GF for the first part $uAd$ of the
decomposition, where the reference line can be anywhere to the
right of the left end of the path (Figure~\ref{fig:prod}).
Similarly, let \bea H_2(x,q,t,v):=\us{r\geq
-n}{\sum_{n\geq0}}\sum_{B\in\D_n}x^{\mathrm{ct}_r(B)}q^{\mathrm{lt}_r(B)}v^r
t^n\label{h2def}\eea be the GF for the second part $B$ of the
decomposition, where now the reference line can be anywhere to the
left of the right end of the path.

\begin{figure}[hbt]
\epsfig{file=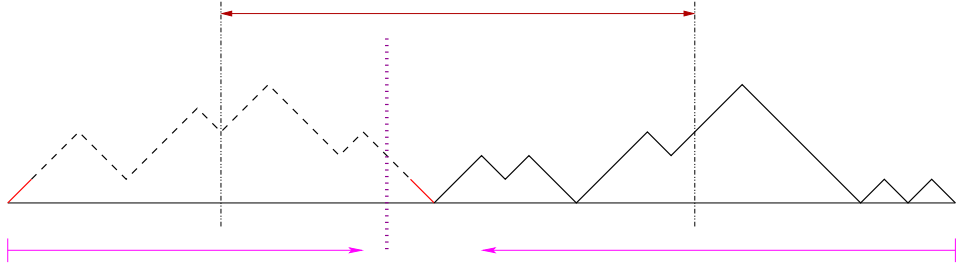,height=3.5cm}
\caption{\label{fig:join} Terms with equal exponent in $t$ and $v$.}
\end{figure}

We would like to express the generating function for paths $uAdB$
in terms of $H_1$ and $H_2$. The product of these two GFs counts
pairs $(uAd,B)$, but if we want the reference line to coincide in
$uAd$ and in $B$, then the two parts are not necessarily placed
next to each other (Figure~\ref{fig:prod}). The exponent of $v$ in
$H_1$ indicates how far to the right the reference line is from
the middle of the path $uAd$. The exponent of $v$ in $H_2$
indicates how far to the left the reference line is from the
middle of the path $B$. In the product $H_1 H_2$, the exponent of
$v$ is the distance from the middle of the path $uAd$ to the
middle of the path $B$ if we draw them so that the reference lines
coincide. Now comes one of the key points of the argument. The
terms that correspond to an actual path $D=uAdB$ are those in
which the two parts are placed next to each other in the picture
($B$ begins  where $uAd$ ends), and this happens precisely when
the exponent of $v$ is half the sum of lengths of $uAd$ and $B$.
But the semilength of each path is the exponent of $t$ in the
corresponding GF, so the sum of semilengths is the exponent of $t$
in the product $H_1 H_2$. Hence, the terms that correspond to
actual paths $D=uAdB$ are exactly those in which the exponent of
$v$ equals the exponent of $t$ (Figure~\ref{fig:join}). In
generating function terminology, the GF consisting of only such
terms is called a \emph{diagonal}.

We also need another variable $y$ to mark the distance between the
reference line and the middle of the new path $D=uAdB$.
Considering that $D$ starts at $(0,0)$, the $x$-coordinate of the
middle of the new path is given by the exponent of $t$ in the
product, which is the sum of the exponents of $t$ in $H_1$ and
$H_2$. The $x$-coordinate of the reference line is given by the
exponent of $t$ in $H_1$ plus the exponent of $v$ in $H_1$. Hence,
the difference between these two $x$-coordinates is given by the
exponent of $t$ in $H_2$ minus the exponent of $v$ in $H_1$.

Let $$P(x,q,t,v,y):=H_1(x,q,t,\frac{v}{y})H_2(x,q,t y,v),$$ and
let its series expansion in $v$ and $t$ be $$P(x,q,t,v,y)=\us{j\ge
-n}{\sum_{n\geq0}}P_{j,n}(x,q,y)v^j t^n.$$ The diagonal (in $v$
and $t$) of $P$ is defined by $$\mathrm{diag}_{v,t}^z\
P:=\sum_{n\ge0}P_{n,n}(x,q,y)z^n.$$ Now, the above argument
implies that this diagonal equals precisely \bea
H_3(x,q,z,y):=\us{-n\le r\le
n}{\sum_{n\ge1}}\sum_{D\in\D_n}x^{\mathrm{ct}_r(D)}q^{\mathrm{lt}_r(D)}y^r
z^n,\label{h3def}\eea that is, the sum over arbitrary non-empty
(since $uAd$ was non-empty) Dyck paths $D$, where the reference
line can be anywhere between the left end and the right end of the
path.

We have found an equation that relates $H_1$, $H_2$ and $H_3$,
thus proving the following lemma.

\begin{lemma} Let $H_1$, $H_2$ and $H_3$ be defined respectively
by (\ref{h1def}), (\ref{h2def}), and (\ref{h3def}). Then,
\bea\mathrm{diag}_{v,t}^z\ H_1(x,q,t,\frac{v}{y})H_2(x,q,t
y,v)=H_3(x,q,z,y).\label{cond}\eea
\end{lemma}

\ms

The next step is to express these three GFs in terms of $G$, so
that (\ref{cond}) will in fact give an equation for $G$.
First, note that given $D\in\D_n$, if $D^R$ is the Dyck path
obtained by reflecting $D$ onto the vertical line $x=n$, then we have that $\mathrm{ct}_{-r}(D)=\mathrm{ct}_r(D^R)$ and
$\mathrm{lt}_{-r}(D)=n-\mathrm{lt}_r(D^R)-\mathrm{ct}_r(D^R)$. Thus,
\bea \sum_{n,r\ge0}\sum_{D\in\D_n}x^{\mathrm{ct}_{-r}(D)}q^{\mathrm{lt}_{-r}(D)}v^r
t^n
=\sum_{n,r\ge0}\sum_{D\in\D_n}\left(\frac{x}{q}\right)^{\mathrm{ct}_r(D^R)}\left(\frac{1}{q}\right)^{\mathrm{lt}_r(D^R)}v^r
(qt)^n=G(\frac{x}{q},\frac{1}{q},qt,v).\hspace*{-12mm}\label{dual}\eea
Also, note that if $|D|=n$ and $r\ge n$, then $\mathrm{ct}_r(D)=\mathrm{lt}_r(D)=\mathrm{ct}_{-r}(D)=0$ and
$\mathrm{lt}_{-r}(D)=n$. In particular,
\bea\us{r>n}{\sum_{n\ge0}}\sum_{D\in\D_n}x^{\mathrm{ct}_r(D)}q^{\mathrm{lt}_r(D)}v^r
t^n=\us{r>n}{\sum_{n\ge0}}C_n v^r
t^n=\sum_{n\ge0}C_n\frac{v^{n+1}}{1-v}t^n=\frac{v}{1-v}C(tv).\label{trivial1}\eea

For $H_1$ we can write \bea
H_1(x,q,t,v)=\us{k\ge-n-1}{\sum_{n\ge0}}\sum_{A\in\D_{n}}x^{\mathrm{ct}_{-k}(uAd)}q^{\mathrm{lt}_{-k}(uAd)}v^k
t^{n+1}=t\left[\us{k>0}{\sum_{n\ge0}}\sum_{A\in\D_{n}}x^{\mathrm{ct}_{-k}(uAd)}q^{\mathrm{lt}_{-k}(uAd)}v^k
t^n\right.\nn\\ \left.
+\sum_{n\ge0}\sum_{A\in\D_{n}}x^{\mathrm{ct}_{0}(uAd)}q^{\mathrm{lt}_{0}(uAd)}t^n
+\us{0<r\le
n+1}{\sum_{n\ge0}}\sum_{A\in\D_{n}}x^{\mathrm{ct}_{r}(uAd)}q^{\mathrm{lt}_{r}(uAd)}v^{-r}
t^n\right]\label{sum3}.\eea For $k>0$, $\mathrm{ct}_{-k}(uAd)=\mathrm{ct}_{-k}(A)$
and $\mathrm{lt}_{-k}(uAd)=\mathrm{lt}_{-k}(A)+1$, so the first sum on the right hand side of (\ref{sum3})
equals \bea
q\us{k>0}{\sum_{n\ge0}}\sum_{A\in\D_{n}}x^{\mathrm{ct}_{-k}(A)}q^{\mathrm{lt}_{-k}(A)}v^k
t^n=q\left[\us{k\ge0}{\sum_{n\ge0}}\sum_{A\in\D_{n}}x^{\mathrm{ct}_{-k}(A)}q^{\mathrm{lt}_{-k}(A)}v^k
t^n\right.\nn\\
\left.-\sum_{n\ge0}\sum_{A\in\D_{n}}x^{\mathrm{ct}_{0}(A)}q^{\mathrm{lt}_{0}(A)}t^n\right]=
q\left[G(\frac{x}{q},\frac{1}{q},qt,v)-G(x,q,t,0)\right],\nn\eea
by (\ref{dual}). For the second sum in (\ref{sum3}), note that
$\mathrm{ct}_{0}(uAd)=\mathrm{ct}_{0}(A)+1$ and $\mathrm{lt}_{0}(uAd)=\mathrm{lt}_{0}(A)$, so the sum
equals
$$x\sum_{n\ge0}\sum_{A\in\D_{n}}x^{\mathrm{ct}_{0}(A)}q^{\mathrm{lt}_{0}(A)}t^n=x G(x,q,t,0).$$
Using that for $r>0$ $\mathrm{ct}_r(uAd)=\mathrm{ct}_r(A)$ and $\mathrm{lt}_r(uAd)=\mathrm{lt}_r(A)$,
the third sum in (\ref{sum3}) can be written as \bea
\us{r>0}{\sum_{n\ge0}}\sum_{A\in\D_{n}}x^{\mathrm{ct}_{r}(A)}q^{\mathrm{lt}_{r}(A)}v^{-r}
t^n-\us{r>n+1}{\sum_{n\ge0}}\sum_{A\in\D_{n}}x^{\mathrm{ct}_{r}(A)}q^{\mathrm{lt}_{r}(A)}v^{-r}
t^n\nn\\
=G(x,q,t,v^{-1})-G(x,q,t,0)-\frac{1}{v(v-1)}C(tv^{-1}),\nn\eea
by (\ref{trivial1}).
Thus, \bea H_1(x,q,t,v)=t\left[q
G(\frac{x}{q},\frac{1}{q},qt,v)+(x-q-1)G(x,q,t,0)+G(x,q,t,\frac{1}{v})
+\frac{1}{v(1-v)}C(\frac{t}{v})\right].\hspace*{-12mm}\label{h1}\eea

For $H_2$, a very similar reasoning implies that
\bea
H_2(x,q,t,v)=G(x,q,t,v)-G(x,q,t,0)+G(\frac{x}{q},\frac{1}{q},qt,\frac{1}{v})+\frac{1}{1-v}C(\frac{qt}{v}).
\label{h2}\eea

Finally, for $H_3$ we get that \bea\label{h3} H_3(x,q,z,y)=
G(x,q,z,y)+G(\frac{x}{q},\frac{1}{q},qz,\frac{1}{y})-G(x,q,z,0)-\frac{y}{1-y}C(zy)
+\frac{1}{1-y}C(\frac{qz}{y})-1. \hspace*{-14mm}\eea

Substituting these expressions for $H_1$, $H_2$ and $H_3$ in
(\ref{cond}) we obtain an equation for $G$. Note that the common
factor $t$ in $H_1(x,q,t,v)$ guarantees that this equation will
express the coefficients of the series expansion in $z$ of
$H_3(x,q,z,y)$ in terms of coefficients of $G$ of smaller order in the
series expansion in $t$ of $H_1(x,q,t,v)H_2(x,q,t,v)$, so it
uniquely determines $G$ as a GF. The final step of the proof is to
guess an expression for $G$ and check that it satisfies this
equation.

\begin{prop}\label{prop} We have
\bea G(x,q,t,v)=\frac{\displaystyle\frac{1-v+(q-1)tv
C(tv)}{1-v+(q-1)tv F_{321}(1,q,t)}-(x-1)tv
C(tv)}{[1-qt(F_{321}(1,q,t)-1)-xt](1-v)} \label{conj}.\eea
\end{prop}

Before proving this proposition, we observe that it implies
theorem~\ref{th:main}. Indeed, we have by definition
$$G(x,q,t,0)=\sum_{n\geq0}\sum_{D\in\D_n}x^{\mathrm{ct}_0(D)}q^{\mathrm{lt}_0(D)}t^n=F_{132}(x,q,t).$$
But if \ref{prop} holds, then
$$G(x,q,t,0)=\frac{1}{1-qt(F_{321}(1,q,t)-1)-xt}=F_{321}(x,q,t),$$
where the last equality follows from (\ref{F321x1}). So, all that
remains is to prove proposition~\ref{prop}.

\begin{proof}

The computations that follow have been done using {\it Maple}.
Let $\wt{H}_1$,
$\wt{H}_2$ and $\wt{H}_3$ be the expressions obtained respectively
from (\ref{h1}), (\ref{h2}) and (\ref{h3}) when $G$ is substituted
with the expression given in (\ref{conj}). All we have to check is
that $$\mathrm{diag}_{v,t}^z\
\wt{H}_1(x,q,t,\frac{v}{y})\wt{H}_2(x,q,t
y,v)=\wt{H}_3(x,q,z,y).$$ Let
$\wt{P}(x,q,t,v,y):=\wt{H}_1(x,q,t,\frac{v}{y})\wt{H}_2(x,q,t
y,v)$. We want to compute diag$_{v,t}^z\ \wt{P}$. In \cite[chapter
6]{EC2}, a general method is described for obtaining diagonals of
rational functions. This theory does not apply to our function
$\wt{P}$, because it is not rational. However, we will show that
in this particular case we can modify the technique to obtain
diag$_{v,t}^z\ \wt{P}$.

The series expansion of $\wt{P}$ in $v$ and $t$,
$$\wt{P}(x,q,t,v,y)=\us{j\ge -n}{\sum_{n\ge0}}\wt{P}_{j,n}(x,q,y)v^j
t^n=\sum_{n,i\ge0}\wt{P}_{i-n,n}(x,q,y)v^i
\left(\frac{t}{v}\right)^n,$$ converges for $|v|<\beta$,
$|\frac{t}{v}|<\alpha$, if $\alpha,\beta>0$ are taken sufficiently
small. Similarly, $$\mathrm{diag}_{v,t}^z\
\wt{P}=\sum_{n\ge0}\wt{P}_{n,n}(x,q,y)z^n$$ converges for $|z|$
sufficiently small. Fix such a small z with $|z|<\alpha\beta^2$.
The series
$$\wt{P}(x,q,t,\frac{z}{t},y)=\us{j\ge -n}{\sum_{n\ge0}}\wt{P}_{j,n}(x,q,y)z^j
t^{n-j}$$ will converge for $|\frac{z}{t}|<\beta$ and
$|\frac{t^2}{z}|<\alpha$. Regarded as a function of $t$, it will
converge for $|t|$ in the annulus
$\frac{|z|}{\beta}<|t|<\sqrt{\alpha |z|}$, which is non-empty
because $|z|<\alpha\beta^2$. In particular, it converges on some
circle $|t|=\rho$ in the annulus. By \cite[Theorem 1]{HauKla71},
$$\mathrm{diag}_{v,t}^z\ \wt{P}=\frac{1}{2\pi
i}\int_{|t|=\rho}{\wt{P}(x,q,t,\frac{z}{t},y)\frac{dt}{t}}.$$ It
can be checked that the singularities of
$\wt{P}(x,q,t,\frac{z}{t},y)/t$ (as a function of $t$) that lie
inside the circle $|t|=\rho$ are all simple poles. These poles are
$$t_1=0,\ t_2=z,\ t_3=\frac{z}{y},\
t_{4,5}=\frac{(1+q)y\pm(1-q)\sqrt{y(y-4qz)}}{2y(y+z(1-q)^2)}z,\
t_{6,7}=\frac{1+q\pm (1-q)\sqrt{1-4zy}}{2(q+zy(1-q)^2)}z.$$ There
are also branch points for $t=\pm\frac{1}{2}\sqrt{\frac{z}{y}}$
and $t=\pm\frac{1}{2}\sqrt{\frac{z}{qy}}$, but they lie outside
the circle for an appropriate choice of $\rho$ in the annulus
$\frac{|z|}{\beta}<\rho<\sqrt{\alpha |z|}$. The remaining
singularities do not depend on $z$ and lie outside the circle.

So, by the residue theorem, the integral can be obtained by
summing up the residues at the poles inside $|t|=\rho$. Computing
them in {\it Maple}, we see that all the residues are 0 except for
those in $t_2$ and $t_3$. Thus,
$$\mathrm{diag}_{v,t}^z\ \wt{P}=\mathrm{Res}_{t=z\ }\wt{P}(x,q,t,\frac{z}{t},y)\frac{1}{t}+
\mathrm{Res}_{t=\frac{z}{y}\
}\wt{P}(x,q,t,\frac{z}{t},y)\frac{1}{t},$$ and this turns out to
be precisely $\wt{H}_3(x,q,z,y)$.

\end{proof}

\section{Some other bijections involving $\S_n(321)$ and $\D_n$}\label{sec:more_bij}

\begin{figure}[hbt]
\epsfig{file=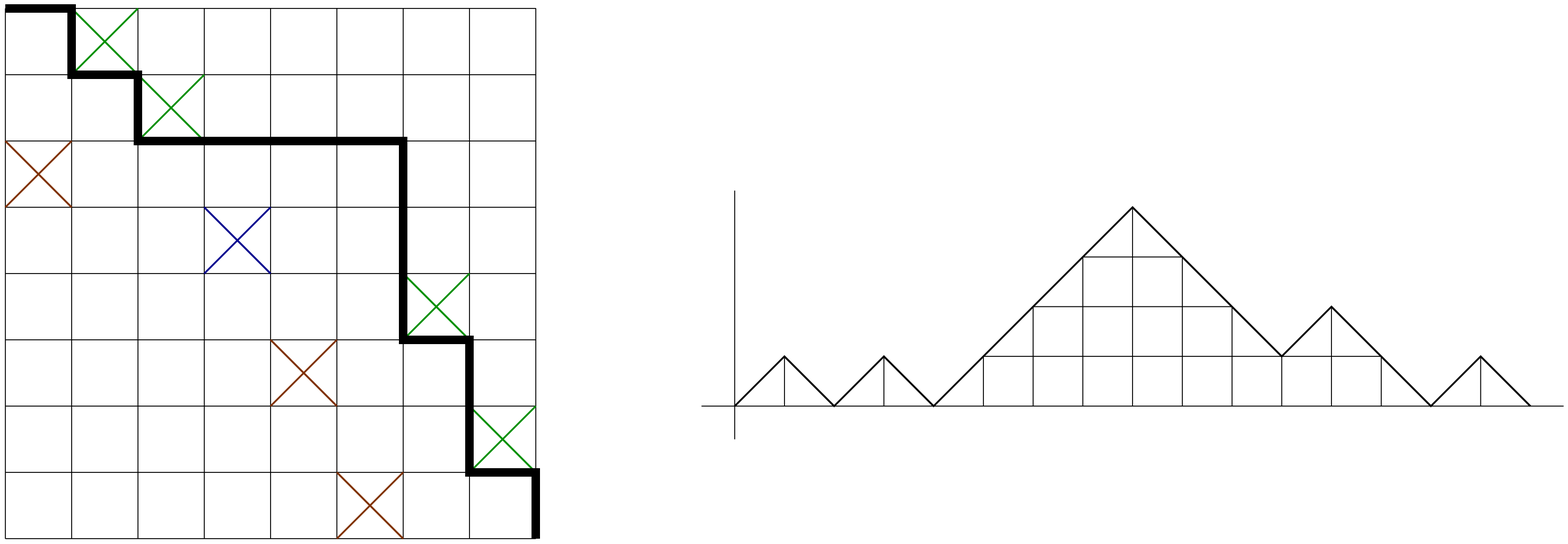,height=3.7cm}
\caption{\label{fig:bij_bjs} The bijection $\bjs$.}
\end{figure}

Looking at permutations as arrays of crosses, as we did to define $\rs$, some other
known bijections between $\S_n(321)$ and $\D_n$ can easily be viewed in a systematic
way, as paths with down and right steps from the upper-left
corner to the lower-right corner of the $n\times n$ array. One such bijection
was established by Billey, Jockusch and Stanley in
\cite[p. 361]{BJS93}. Denote it by $\bjs$. Consider the path that
leaves the crosses corresponding to excedances to the right, and
stays always as far from the main diagonal as possible
(Figure~\ref{fig:bij_bjs}). Then $\bjs(\pi)$ can be obtained from
it just by reading an up-step every time the path moves to the
right and a down-step every time the path moves down.

\begin{figure}[hbt]
\epsfig{file=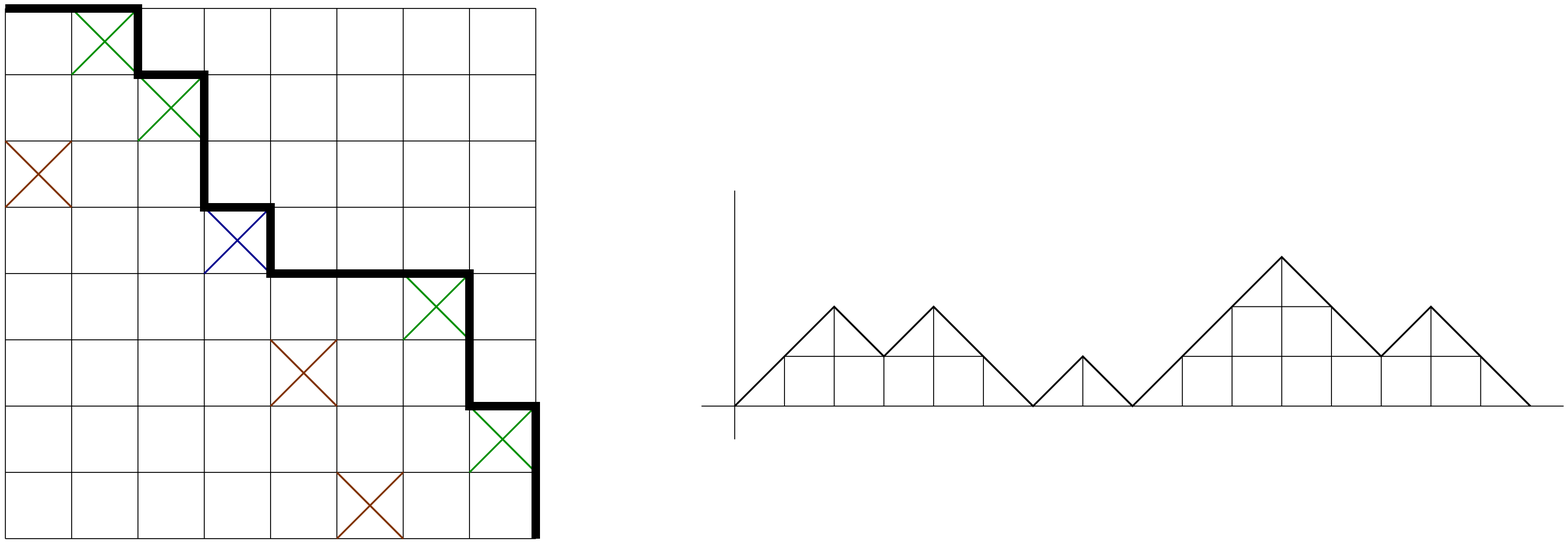,height=3.7cm}
\caption{\label{fig:bij_kra} The bijection $\kra$.}
\end{figure}

In \cite{Kra01}, Krattenthaler describes a bijection from
$\S_n(123)$ to $\D_n$. If we omit the last step, consisting in
reflecting the path into a vertical line, and compose the
bijection with the reversal operation, that maps a permutation
$\pi_1\pi_2\cdots\pi_n$ into $\pi_n\cdots\pi_2\pi_1$, we get a
bijection from $\S_n(321)$ to $\D_n$. Denote it by $\kra$. In the
array representation, $\kra(\pi)$ corresponds (by the same trivial
transformation as before) to the path that leaves all the crosses
to the left and remains always as close to the main diagonal as
possible (Figure~\ref{fig:bij_kra}).

\begin{figure}[hbt]
\epsfig{file=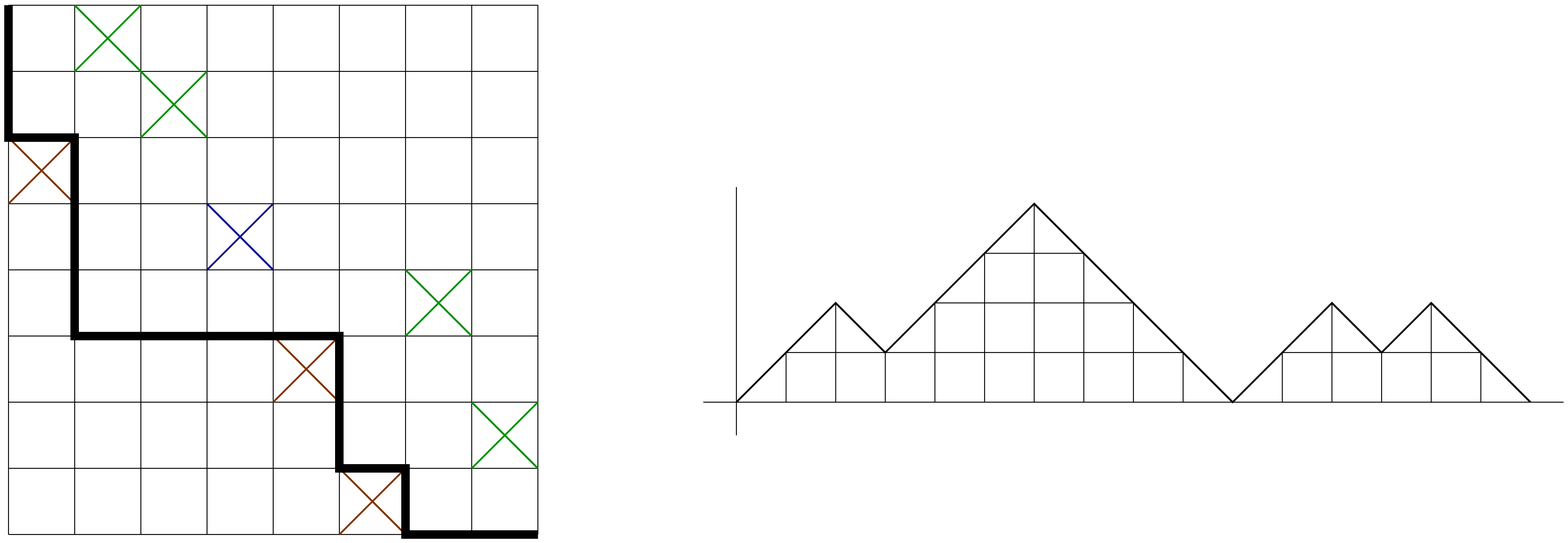,height=3.7cm} \caption{\label{fig:bij_4}
The bijection $\bijn$.}
\end{figure}

Our first bijection is related to this last one by
$\rs(\pi)=\kra(\pi^{-1})$. In a similar way, we could still define
a fourth bijection $\bijn:\S_n(321)\longrightarrow\D_n$ by
$\bijn(\pi):=\bjs(\pi^{-1})$ (Figure~\ref{fig:bij_4}).

Combining these bijections and their inverses, one can get some
automorphisms on Dyck paths and on 321-avoiding permutations with
interesting properties. Recall that a \emph{valley} of a Dyck path
$D$ is a down-step followed by an up-step ($du$ in the Dyck word).
Denote by $\mathrm{va}(D)$ the number of valleys of $D$. Denote by $p_2(D)$
the number of peaks of $D$ of height at least 2. Clearly, both
$p_2(D)+h(D)$ and $\mathrm{va}(D)+1$ equal the total number of peaks of $D$.

It can be checked that $\rs\circ\bjs^{-1}$ is an involution on
$\D_n$ with the property that
$\mathrm{va}(\rs\circ\bjs^{-1}(D))=\mathrm{dr}(D)$ and
$\mathrm{dr}(\rs\circ\bjs^{-1}(D))=\mathrm{va}(D)$. Indeed, this
follows from the fact that excedances are sent to valleys by
$\bjs$ and to double rises by $\rs$. This bijection gives yet
another proof of the symmetry of the bivariate distribution of the
pair $(\mathrm{va},\mathrm{dr})$ of statistics in Dyck paths. A
different involution with this property was introduced in
\cite{Deu99}.

Another involution on $\D_n$ is given by $\rs\circ\kra^{-1}$. This
one shows the symmetry of the distribution of the pair
$(\mathrm{dr},p_2)$, because
$\mathrm{dr}(\rs\circ\kra^{-1}(D))=p_2(D)$ and
$p_2(\rs\circ\kra^{-1}(D))=\mathrm{dr}(D)$. Besides, it preserves
the number of hills, i.e., $h(\rs\circ\kra^{-1}(D))=h(D)$. To see
this, just note that both $\kra$ and $\rs$ send fixed points to
hills, whereas excedances are sent to peaks of height at least 2
by $\kra$ and to double rises by $\rs$.

On the other hand, the involution on $\S_n(321)$ that maps $\pi$
to $(\bjs^{-1}(\rs(\pi)))^{-1}$ gives a combinatorial proof of the
fact that the number of 321-avoiding permutations with $k$
excedances equals the number of 321-avoiding permutations with
with $k+1$ weak excedances (recall that $i$ is a \emph{weak
excedance} of $\pi$ if $\pi_i\ge i$). The analogous result for
general permutations is well known. An implication of
Theorem~\ref{th:main} is that this result is also true for
132-avoiding permutations.

\small
\subsection*{Acknowledgements}

The author is grateful to Richard Stanley for suggesting the
problem that led to Theorem~\ref{th:main}, and for the main
results in section~\ref{sec:321}.

\subsection*{Note}

A bijective proof of Theorem~\ref{th:main} has very recently been found by 
the author and Igor Pak. This result is intended to appear in~\cite{EliPak},
which is currently in preparation.


\begin{thebibliography}{99}

\bibitem{BJS93} S. Billey, W. Jockusch, R. Stanley, Some
Combinatorial Properties of Schubert Polynomials, {\it J. Alg.
Comb.} 2 (1993), 345--374.

\bibitem{Deu99} E. Deutsch, An involution on Dyck paths and its
consequences, {\it Discrete Math.} 204 (1999), 163--166.

\bibitem{Deu99_} E. Deutsch, Dyck Path Enumeration, {\it Discrete Math.}
204 (1999), 167--202.

\bibitem{EliPak} S. Elizalde, I. Pak, Bijections for Refined Restricted 
Permutations (in preparation).

\bibitem{FlSe98}
 P. Flajolet, P. and R. Sedgewick, {\it Analytic combinatorics}
(book in preparation) (1998). (Individual chapters are available as
INRIA Research Reports 1888,
  2026, 2376, 2956, 3162.).

\bibitem{HauKla71} M.L.J. Hautus, D.A. Klarner, The diagonal of a
double power series, {\it Duke Math. J.} 38, No.2 (1971).

\bibitem{Knu73} D. Knuth, The Art of Computer Programming, Vol. I
(Addison-Wesley, Reading, MA, 2nd ed., 1973).

\bibitem{Kra01} C. Krattenthaler, Permutations with restricted
patterns and Dyck paths, {\it Adv. Appl. Math.} 27 (2001),
510--530.

\bibitem{Rei02} A. Reifegerste, On the diagram of 132-avoiding
permutations, arxiv:math.CO/0208006 v3  15 Oct 2002.

\bibitem{Ric88} D. Richards, Ballot sequences and restricted
permutations, {\it Ars Combin.} 25 (1988), 83--86.

\bibitem{RSZ02} A. Robertson, D. Saracino, D. Zeilberger, Refined
Restricted Permutations, arxiv:math.CO/0203033 v1  4 Mar 2002.

\bibitem{SeFl96}  R. Sedgewick and P. Flajolet,
{\it An introduction to the analysis of algorithms}
(Addison-Wesley, 1996).

\bibitem{SS85} R. Simion, F.W. Schmidt, Restricted Permutations,
{\it European J. Combin.} 6 (1985), 383--406.

\bibitem{EC1} R. Stanley, {\it Enumerative Combinatorics}, vol. I
(Cambridge Univ. Press, Cambridge, 1997).

\bibitem{EC2} R. Stanley, {\it Enumerative Combinatorics}, vol. II
(Cambridge Univ. Press, Cambridge, 1999).

\bibitem{Wes95} J. West, Generating trees and the Catalan and Schr\"oder numbers,
{\it Discrete Math.} 146 (1995), 247--262

\bibitem{Wes96} J. West, Generating Trees and Forbidden
Subsequences, {\it Discrete Math.} 157 (1996), 363--374.

\end{thebibliography}
\end{document}